\documentclass{article}

\usepackage[utf8]{inputenc}
\usepackage[T1]{fontenc}
\usepackage{geometry,amsmath,amssymb,amscd}

\usepackage{enumerate}

\newtheorem{thm}{Theorem}[section]

\newtheorem{prop}[thm]{Proposition}

\newcommand \cL {\mathcal {L}}

\newcommand \RR{\mathbb{R}}

\DeclareMathOperator{\bb}{\quad \Box}

\begin{document}

\title {A criterion for the differential  flatness of a nonlinear control system}

\author {Bruno Sauvalle}

\author {Bruno Sauvalle \\ MINES ParisTech, PSL - Research University \\60 Bd Saint-Michel, 75006 Paris, France }
\date{November 2017}

\maketitle

\AtEndDocument{\bigskip{\footnotesize%
  \textsc{Bruno Sauvalle, MINES ParisTech, PSL - Research University, 60 Bd Saint-Michel, 75006 Paris, France } \par  
  \textit{E-mail address} \texttt{bruno.sauvalle@mines-paristech.fr}  
}}

\maketitle
\begin{abstract}
Let's consider a control system described by the implicit equation $F(x,\dot x) = 0$. If this system is differentially flat, then the following criterion is satisfied : For some integer $r$, there exists a function $\varphi(y_0, y_1, ..,y_r)$ satisfying the following conditions: (1) The map $(y_0,..,y_{r+1}) \mapsto ( \varphi( y_0, y_1, ..,y_r), \frac {\partial \varphi}{\partial y_0}y_1 +\frac {\partial \varphi}{\partial y_1}y_2+ ..+ \frac {\partial \varphi}{\partial y_r} y_{r+1})$ is a submersion on the variety $F(x,p) = 0$. (2) The map $y_0 \mapsto x_0 = \varphi(y_0,0,..,0)$ is a diffeomorphism on the equilibrium variety $F(x,0) = 0$.

 Inversely, if a control system satifies this flatness criterion, then it is  locally controllable  at equilibrium points.
\end{abstract}

\section {Introduction}

	The purpose of this note is to propose a new approach to the study of the differential flatness of a control system. For an overview of flatness theory, see for example $\cite {Levine}$. 
	For an overview of existing differential flatness criteria, see \cite{Martin}, section 3.1.2, and $\cite{Levine2}$.

We  consider a control system defined by the implicit equation $F(x,\dot x) = 0$, where $F(x,p)$ is a smooth function defined on $\RR^n \times \RR^n$, with rank $\frac {\partial F}{ \partial p} = n-m$ ($m>0$). We assume that for each $x$ there exists at least one $p$ satisfying the equation, and more generally that the manifold $F(x,p)$ can be parameterized using a smooth function $u \mapsto f(x,u)$ satisfying $F(x,f(x,u)) = 0$. This condition is trivially satisfied if the system $F(x, \dot x) = 0$ is derived from some  explicit description of the control system such as  $\dot x = f(x,u)$ .

\section {Finding parameter functions }

 We say that a function $\varphi(y_0,..y_r)$ defined on $(\RR^m)^{r+1}$ with values in $\RR^n$ is a parameter function for the system $F$ if  we have for all smooth functions $y(t)$ defined on $\RR$ with values in  $\RR^m$, and writing $x(t) = \varphi ( y(t), \dot y(t), ..,y^{(r)}(t))$, the equality $F( x(t), \frac {d}{dt}x(t)) = 0 $.
 We can also write this condition as : 

\begin {equation} \label {timeq} F( \varphi( y(t), \dot y(t) ..y^{(r)}(t)), \frac {\partial \varphi}{\partial y_0}\dot y(t) +..+ \frac {\partial \varphi}{\partial y_r} y^{(r+1)}(t)) = 0 \end {equation}

Considering that for $t$ fixed, the function $y$ can be chosen so that $y(t), \dot y(t), ..,y^{(r+1)}(t)$ take any independent values, we see that we have for any vectors $y_0, .., y_{r+1}$ in $\RR^m$ the equality 
\begin {equation}  \label {PDE} F( \varphi( y_0, y_1, ..,y_r), \frac {\partial \varphi}{\partial y_0}y_1 +\frac {\partial \varphi}{\partial y_1}y_2 ..+ \frac {\partial \varphi}{\partial y_r} y_{r+1}) = 0 \end {equation}

Inversely, it is clear that if a function $\varphi$ satisfies this partial differential equation, then equation ($\ref {timeq}$) is satisfied for all smooth functions $y$. 

The question of the existence and computation of parameter functions is then equivalent to the question of existence and computation of solutions to this implicit first order partial differential equation. A wide variety of techniques are available to study this question : Cartan K\"ahler theory, method of characteristics and  Charpit equations if $F$ is scalar, etc ..

The main difficulty here is that $\varphi$ has to be a function of $y_0,..y_r$ only and must not depend of $y_{r+1}$, which appears as a free parameter in the equation. This requires that the ruled manifold criterion described in $\cite{Rouchon}$ be satisfied.

\section { Differentially flat systems}

Let's recall (cf \cite {Levine}) that a system $F(x,\dot x)$ is called differentially flat if there exists a parameter function $\varphi$ satisfying equation $(\ref{timeq})$ and another function $\psi(x_0,..,x_s)$ so that a solution $x(t)$ of the system can be written as $x(t) = \varphi( y, \dot y ... y^{(r)})$ if and only if  $y(t) = \psi( x, \dot x ,.. x^{(s)})$. This requirement is  stronger than the existence of a parameter function. Let's note $\cL_\tau$ for the derivation operator $\cL_\tau \varphi (y) = \sum_{i = 0}^r\frac {\partial \varphi }{\partial y_i}y_{i+1}$.
If $\varphi$ is associated to a flat system, we have :
\begin {prop}\label {surj1}
The smooth map $\Phi : (y_0,..,y_{r+1}) \mapsto (\varphi(y_0..y_r), \cL_\tau \varphi(y_0..y_{r+1}))$  is a submersion on the manifold $F(x,p) = 0$.  \end {prop}
Remark : this implies that $d\varphi$ is surjective.\\
Proof : We have to prove that $\Phi$ and $d\Phi$ are surjective.

Let's consider an element $(x_0, p_0)$ of the manifold and the associated parameter $u_0$ so that $p_0 = f(x_0,u_0)$. Let's now choose the constant function $u(t)$ defined by the formula  $u(t) = u_0$. Then the unique solution of the firts order ODE $\dot x = f(x,u)$ with $\dot x(0) = p_0$, $x(0) = x_0$ satisfies $F(x,\dot x) = 0$ for all values of $t$ and is a smooth function of $x_0$ and $p_0$. Using the function $\psi$, we get some function $y(t)$ which is also a smooth function of $x_0$ and $p_0$ so that $x(t) = \varphi(y, \dot y,..y^{(r)})$.
Writing $y(0) = y_0, ..,y^{(r+1)}(0) = y_{r+1}$ can then write that 
$$ (x_0,p_0) = (\varphi(y_0,..,y_r),\cL_\tau \varphi(y_0..y_{r+1}) = \Phi(y_0,..y_{r+1}) $$
This shows that $\Phi$ is surjective.
Taking the differential of this equation with respect to $(x_0,p_0)$, we get 
$$ Id = d\Phi \circ (..)$$
which shows that $d\Phi$ is surjective $\bb$

\begin {prop} \label {eq1}
The restriction of $\Phi$ to the elements of the form $(y_0,0..0)$ is a diffeomorphism on the equilibrium manifold $F(x_0,0) = 0$
\end {prop}
Proof : It is well known that Lie-B\"acklund equivalence preserves equilibrium points ( cf \cite{Levine}, Th 5.2) $\bb$

\section {Controllability at equilibrium points }

We now consider in this section a control system described by an equation $F(x,\dot x) = 0$ such that the equation $ (\ref{PDE})$ has a solution $\varphi$ satisfying propositions \ref{surj1}  and \ref{eq1}.

\begin {prop}
 The system is locally controllable at equilibrium points. 
\end {prop}

Proof : We consider some equilibrium point $x_0$, so that $F(x_0,0) = 0$. Using $\ref{surj1}$, we see that $(x_0,0)$ is in the image of $\Phi$ and that $d\varphi$ is surjective at this point.

We now show that Kalman criterion is satisfied. Let's note $A =  \frac {\partial f}{\partial x} $ and $B = Im  \frac {\partial f}{\partial u}$
Let's take the derivative of the equation $F(x,f(x,u)) = 0$
We get 
$$ \frac {\partial F}{\partial x} + \frac {\partial F}{\partial p} A = 0 $$
and, considering that for $x$ fixed, the map $u \mapsto f(x,u)$ is assumed to be a diffeomorphism on the manifold $F(x,p) = 0$,
$$ B = Ker \frac {\partial F}{\partial p} $$

We  can  write the equation (\ref{PDE}) as   

$$ F( \varphi, \cL_\tau \varphi) = 0$$
Taking differentials with respect to $y_0..y_{r+1}$, we get the equalities : \\

$$ \frac {\partial F}{\partial x} \frac {\partial \varphi}{\partial y_i} + \frac {\partial F}{\partial p} \frac {\partial  }{\partial y_i}\cL_\tau \varphi  = 0 $$
using the commutation relation $[ \frac {\partial} {\partial y_i}, \cL_\tau] = \frac {\partial}{\partial y_{i-1}}$ for $1 \le i \le r+1$ and $[ \frac {\partial} {\partial y_0}, \cL_\tau] = 0$, we get \\
\begin {itemize}

\item Differential with respect to $y_0$ : 

$$   \frac {\partial F}{\partial x}  \frac {\partial \varphi}{\partial y_0 }   +  \frac {\partial F}{\partial p} \cL_\tau  \frac {\partial \varphi}{\partial y_0 } = 0 $$
\item Differential with respect to $y_i$ for  $i \in 1..r$ : 
$$  \frac {\partial F}{\partial x}  \frac {\partial \varphi}{\partial y_i }   +  \frac {\partial F}{\partial p} \cL_\tau  \frac {\partial \varphi}{\partial y_i } +  \frac {\partial F}{\partial p}\frac {\partial \varphi}{\partial y_{i-1} } = 0 $$

\item  Differential with respect to  $y_{r+1}$ ( considering that $\varphi$ is a function of $y_0..y_r$ only) :
$$  \frac {\partial F}{\partial p} \frac {\partial \varphi}{\partial y_{r} } = 0 $$

\end {itemize}

Considering that $(x_0,0)$ is an equilibrium point of the system, we now use $\ref{eq1}$ to get  $y_1 = y_2 = ..= y_r= y_{r+1} = 0$. 
The equations become : \\
\begin {itemize}
\item  Differential with respect to $y_0$ : 

$$   \frac {\partial F}{\partial x}  \frac {\partial \varphi}{\partial y_0 }   = 0$$

\item  Differential with respect to $y_i$ for  $i \in 1..r$ : 
$$  \frac {\partial F}{\partial x}  \frac {\partial \varphi}{\partial y_i }   +  \frac {\partial F}{\partial p}\frac {\partial \varphi}{\partial y_{i-1} } = 0 $$

\item Differential with respect to  $y_{r+1}$ : 
$$   \frac {\partial F}{\partial p}\frac {\partial \varphi}{\partial y_{r} } = 0 $$
\end {itemize}

Writing $ \frac {\partial F}{\partial x} = -\frac {\partial F}{\partial p} A$, 
The equations become 

$$   \frac {\partial F}{\partial p}(A  \frac {\partial \varphi}{\partial y_0 } )  = 0 $$

$$  \frac {\partial F}{\partial p}(-A  \frac {\partial \varphi}{\partial y_i }   +  \frac {\partial \varphi}{\partial y_{i-1} }) = 0 $$

$$  \frac {\partial F}{\partial p}\frac {\partial \varphi}{\partial y_{r} } = 0 $$

this implies inductively the inclusions 
$$Im \frac {\partial \varphi}{\partial y_{r} } \subset B $$
$$Im  \frac {\partial \varphi}{\partial y_{r-1} } \subset AB+B $$
..
$$ Im \frac {\partial \varphi}{\partial y_{0}} \subset A^rB+A^{r-1}B+..+B $$

and we see that $A^rB+A^{r-1}B+..+B $ has to contain the range of all the matrices $\frac {\partial \varphi}{\partial y_{i} }$, which generate $\RR^n$ $\bb$


\begin{thebibliography}{HD}

\bibitem {Levine} J. Levine,
\emph{ Analysis and Control of Nonlinear Systems : A Flatness based approach}, 
 Springer, Mathematical Engineering, 2009
 
 \bibitem {Levine2}J. Levine
 \emph {On necessary and sufficient conditions for differential flatness}
 J. AAECC (2011) 22:47
 
 \bibitem {Martin} Ph. Martin, R.M. Murray, and P. Rouchon. 
 \emph{ Flat systems},
  In G. Bastin and M. Gevers, editors,
 Plenary Lectures and Minicourses, Proc ECC 97, Brussels, pages 211-264, 1997
 


\bibitem{Rouchon} P. Rouchon, 
\emph {Necessary Condition and Genericity of Dynamic Feedback Linearization}, 
J. of Math. Systems, Estim. and Control 5(3),1995, pp. 345-358

\end{thebibliography}
\end{document}